\documentclass[11pt]{amsart}

\usepackage[T1]{fontenc}
\usepackage[utf8]{inputenc}
\usepackage{lmodern}
\usepackage[letterpaper,margin=1in]{geometry}
\usepackage{mathtools,amssymb}
\usepackage{microtype}
\usepackage{enumitem}
\usepackage{hyperref}
\hypersetup{
  pdftitle={An affine reverse Faber--Krahn inequality via information decay},
  pdfauthor={Matthew J. Colbrook and George Stepaniants},
  colorlinks=true,
  linkcolor=blue,
  citecolor=blue,
  urlcolor=blue
}
\allowdisplaybreaks[1]
\numberwithin{equation}{section}

\newtheorem{theorem}{Theorem}[section]
\newtheorem{proposition}[theorem]{Proposition}
\newtheorem{lemma}[theorem]{Lemma}

\theoremstyle{definition}

\theoremstyle{remark}
\newtheorem{remark}[theorem]{Remark}

\newcommand{\R}{\mathbb R}
\newcommand{\E}{\mathbb E}
\newcommand{\Prob}{\mathbb P}
\newcommand{\Id}{I_n}
\newcommand{\gam}{\gamma_n}
\newcommand{\CP}{C_{\mathrm P}}
\newcommand{\HS}{\mathrm{HS}}
\newcommand{\op}{\mathrm{op}}
\newcommand{\eps}{\varepsilon}
\newcommand{\norm}[1]{\lVert #1\rVert}
\DeclareMathOperator{\Cov}{Cov}
\DeclareMathOperator{\Var}{Var}
\DeclareMathOperator{\supp}{supp}
\DeclareMathOperator{\Div}{div}
\DeclareMathOperator{\tr}{tr}
\DeclareMathOperator{\SL}{SL}
\DeclareMathOperator{\Int}{int}

\title[Affine reverse Faber--Krahn]{An affine reverse Faber--Krahn inequality\\via information decay}
\author{Matthew J. Colbrook}
\address{Department of Applied Mathematics and Theoretical Physics,
University of Cambridge, Wilberforce Road, Cambridge CB3 0WA, United Kingdom}
\email{m.colbrook@damtp.cam.ac.uk}
\thanks{Matthew J. Colbrook is the corresponding author.}
\author{George Stepaniants}
\address{Department of Applied Mathematics and Theoretical Physics,
University of Cambridge, Wilberforce Road, Cambridge CB3 0WA, United Kingdom}
\curraddr{Courant Institute of Mathematical Sciences, New York University,
251 Mercer Street, New York, NY 10012, USA}
\email{gs852@cam.ac.uk}
\date{27 October 2026}

\begin{document}

\begin{abstract}
We prove that every convex body \(K\subset\R^n\) admits
\(S\in\SL_n(\R)\) such that
\(\lambda_1(SK)|K|^{2/n}\le Cn\), where \(C\) is universal.
For origin-symmetric bodies of volume one, the bound holds in maximal
Gaussian measure position. The main step is the spectral estimate
\[
\lambda_1(D)\le C(P)\bigl(n+\sqrt{n I(\nu\mid\gam)}\bigr)
\]
for an even isotropic log-concave probability measure \(\nu\) supported on
a symmetric convex body \(D\), with Poincar\'e constant at most \(P\) and
finite relative Fisher information.
A second-order Stein estimate and an information jump under convolution
give entropy decay faster than \(m^{-1}\) for normalised sums. Their
supports control the Gaussian measure of dilates, and a Brownian bridge
expansion gives the spectral estimate. A Gaussian-weighted ground state
supplies a measure whose Fisher information is bounded in terms of
\(\lambda_1(K)\). The square-root dependence then gives the bound of
order \(n\).
\end{abstract}

\maketitle

\section{Introduction}

Let \(n\ge1\) be an integer and \(K\subset\R^n\) be a convex body (a
compact convex set with nonempty interior). The first eigenvalue of the Dirichlet Laplacian \(-\Delta\) is
\[
\lambda_1(K)
=\inf_{0\ne f\in H_0^1(\Int K)}
  \frac{\int_{\Int K}|\nabla f|^2\,\mathrm{d}x}
       {\int_{\Int K}f^2\,\mathrm{d}x}.
\]
The affine reverse Faber--Krahn problem asks whether there exists a universal constant $C$ such that
\[
\inf_{S\in\SL_n(\R)}\lambda_1(SK)\,|K|^{2/n}\le Cn.
\]
This question goes back to
Schmuckenschl\"ager~\cite{Sch11}; see also
\cite[Conjecture~5.1]{FP26}. An earlier bound of order \(n\log^2 n\)
is due to Klartag and E.~Milman; see
\cite[Proposition~4.8]{BK26} for the argument and its history.
Faifman and Polterovich~\cite[Theorem~5.2]{FP26} obtain the bound
\(Cn\log(2+n)\).

We prove the conjectured bound of order \(n\). A formal verification in Lean of the theorems of this paper is provided in \cite{CS26}.

\begin{theorem}\label{thm:affine-main}
There is a universal constant \(C>0\) such that every convex body
\(K\subset\R^n\) admits a map \(S\in\SL_n(\R)\) satisfying
\[
\lambda_1(SK)\,|K|^{2/n}\le Cn.
\]
\end{theorem}

For origin-symmetric bodies, the estimate holds in a classical affine
position defined by Gaussian measure.
Let \(\gam\) be standard Gaussian probability measure on \(\R^n\), with
density
\[
\varphi_n(x)=(2\pi)^{-n/2}e^{-|x|^2/2}.
\]
For a body \(K\), let
\[
\gamma_K(A)=\frac{\gam(A\cap K)}{\gam(K)}.
\]
An origin-symmetric body \(K\) of volume one is in
\emph{maximal Gaussian measure position} if
\[
\gam(K)=\max_{S\in\SL_n(\R)}\gam(SK).
\]
Such a position exists. Bobkov~\cite{Bob11} showed that the restricted
Gaussian has scalar covariance bounded above and below by absolute
constants in this position; see Section~\ref{sec:gaussian-position}.
Recent applications of this position to mean width and metric entropy
appear in Paouris and Pathak~\cite{PP26}.

\begin{theorem}\label{thm:gaussian-main}
There is a universal constant \(C>0\) such that every origin-symmetric
convex body \(K\subset\R^n\) with \(|K|=1\), in maximal Gaussian measure
position, satisfies
\[
\lambda_1(K)\le Cn.
\]
\end{theorem}

\subsection{Outline of the proof}

The Stein estimates used in our proof come from the spectral-gap
construction of Courtade, Fathi and Pananjady~\cite{CFP19} and its
higher-order extension by Fathi~\cite{Fathi21}. Their entropy estimates
build on the entropy--Stein--information method of Ledoux, Nourdin and
Peccati~\cite{LNP15}. We combine the second-order estimate with Bizeul's
information jump~\cite{Biz23}, a Fisher-information counterpart of the
entropy-jump inequality of Ball and Nguyen~\cite{BN12}.
Our proof uses the resulting entropy decay to bound
\(-\log\gam(tD)\) with a decay exponent greater than two, and then
computes the survival probability of a Brownian bridge giving the
following spectral estimate.

\begin{theorem}\label{thm:information-spectrum}
Let \(D\subset\R^n\) be an origin-symmetric convex body, and let
\(\nu\) be an even isotropic log-concave probability measure supported
on \(D\). Suppose that
\[
\CP(\nu)\le P,\qquad P\ge2,
\]
and that \(I_0=I(\nu\mid\gam)<\infty\).
Then
\begin{equation}\label{eq:information-spectrum}
\lambda_1(D)\le C(P)\left(n+\sqrt{nI_0}\right),
\end{equation}
where \(C(P)\) is finite and depends only on \(P\).
\end{theorem}

For independent \(U_1,\ldots,U_m\) with law \(\nu\), let
\(\nu_m=\mathcal L(m^{-1/2}\sum_{j=1}^m U_j)\).
The second-order Stein construction in Section~\ref{sec:stein}, combined
with Bizeul's information jump, gives
\[
H(\nu_m\mid\gam)\le C(P)\sqrt{nI_0}\,m^{-1-\delta/2}
\]
for dyadic \(m\), where \(\delta=\delta(P)>0\).
Since \(\nu_m\) is supported on \(\sqrt mD\), this implies
\[
-\log\gam(tD)\le C(P)\sqrt{nI_0}\,t^{-2-\delta},\qquad t\ge1.
\]
The exponent greater than two allows the coefficients of a dyadic
Brownian bridge expansion to satisfy both a summable probability bound
and a summable uniform bound. The resulting survival estimate bounds the
Dirichlet heat kernel from below and proves
Theorem~\ref{thm:information-spectrum}.

To apply the theorem, we need a measure with bounded Poincar\'e constant
and controlled covariance, whose Fisher information may depend on the
eigenvalue. The restricted Gaussian in maximal Gaussian measure position
has the first two properties, but its density has a jump at the boundary.
In Section~\ref{sec:main-proof}, we convolve it with the probability
measure whose density is proportional to \(\psi^2e^{-|x|^2/2}\), where
\(\psi\) is a first Dirichlet eigenfunction. After covariance
normalisation, the resulting measure satisfies
\[
\CP(\nu)\le C,\qquad I(\nu\mid\gam)\le8\lambda_1(K)+n.
\]
Its support is a linear image of \(2K\), with uniformly bounded
distortion. Theorem~\ref{thm:information-spectrum} therefore gives
\[
\lambda_1(K)\le C\left(n+\sqrt{n\lambda_1(K)}\right).
\]
Absorbing the square-root term proves Theorem~\ref{thm:gaussian-main}.
Volume normalisation and the symmetric-intersection argument complete
the proof of Theorem~\ref{thm:affine-main}.

\subsection{Notation and conventions}

A probability measure is \emph{isotropic} if it has mean zero and
covariance \(\Id\). It is \emph{even} if it is invariant under
\(x\mapsto-x\). The law of a random vector \(U\) is denoted by
\(\mathcal L(U)\). Quantities such as \(\CP(U)\) and \(J(U)\) refer
to its law. A dyadic integer is \(m=2^r\) with \(r\in\{0,1,2,\ldots\}\).

For a probability density \(\rho\), define its ordinary Fisher
information by
\begin{equation}\label{eq:ordinary-fisher}
J(\rho\,\mathrm{d}x)=
\begin{cases}
4\displaystyle\int_{\R^n}|\nabla\sqrt{\rho}|^2\,\mathrm{d}x,
 &\sqrt{\rho}\in H^1(\R^n),\\
+\infty,&\text{otherwise}.
\end{cases}
\end{equation}
For smooth positive densities this equals
\(\int|\nabla\log\rho|^2\rho\,\mathrm{d}x\).
For \(\nu=\rho\,\mathrm{d}x\), its relative entropy with respect to \(\gam\) is
\[
H(\nu\mid\gam)=\int_{\R^n}
 \rho\log\frac{\rho}{\varphi_n}\,\mathrm{d}x,
\]
with \(0\log0=0\); for a measure without a density we set
\(H(\nu\mid\gam)=+\infty\). For a density with finite second moment,
we use the
Sobolev definition
\[
I(\nu\mid\gam)=
\int_{\R^n}|2\nabla\sqrt\rho+x\sqrt\rho|^2\,\mathrm{d}x
\]
when \(\sqrt\rho\in H^1(\R^n)\), and \(I(\nu\mid\gam)=+\infty\)
otherwise. For smooth positive densities this is
\(\int|\nabla\log(\rho/\varphi_n)|^2\rho\,\mathrm{d}x\).
If \(J(\nu)<\infty\), integration by parts with cutoffs gives
\(\int x\cdot\nabla\rho\,\mathrm{d}x=-n\), and hence
\[
I(\nu\mid\gam)=J(\nu)+\int|x|^2\,\mathrm{d}\nu-2n.
\]
In particular, for isotropic \(\nu\),
\begin{equation}\label{eq:relative-ordinary}
I(\nu\mid\gam)=J(\nu)-n.
\end{equation}
Here integration by parts is justified by inserting standard smooth cutoffs at radius \(R\) and letting \(R\to\infty\).
Indeed, Cauchy–Schwarz gives \(x\cdot\nabla\rho\in L^1(dx)\) while the terms involving cutoff derivatives are bounded by \(C\int_{|x|\ge R}\rho,dx\to0\).

The Poincar\'e constant \(\CP(\nu)\) is the least constant such that,
for every smooth function \(f\),
\begin{equation}\label{eq:poincare}
\Var_\nu f\le \CP(\nu)\int|\nabla f|^2\,\mathrm{d}\nu.
\end{equation}
The inequality extends to Sobolev limits. In particular, isotropy implies
\(\CP(\nu)\ge1\).
For symmetric matrices, \(A\preceq B\) means that \(B-A\) is positive
semidefinite. We denote the Euclidean operator norm by
\(\|\cdot\|_{\op}\), and the Hilbert--Schmidt norm on matrices and
tensors by \(\|\cdot\|_{\HS}\). Thus
\(\|T\|_{L^2(\nu)}^2=\int\|T\|_{\HS}^2\,\mathrm{d}\nu\).
The colon denotes coefficientwise contraction. For a vector field
\(F\), our derivative convention is \((DF)_{ij}=\partial_jF_i\).
For a third-order tensor \(T\),
\[
T:D^2F=\sum_{i,j,k}T_{ijk}\,\partial_j\partial_k F_i.
\]
For a matrix \(M\), the vector \(T:M\) has components
\((T:M)_i=\sum_{j,k}T_{ijk}M_{jk}\).

Symbols \(c,C>0\) denote absolute constants whose values may change from
one occurrence to the next. We also denote dependence on a parameter by \(C(P)\), or \(C(p)\).

\section{Gaussian position and elementary analytic facts}
\label{sec:gaussian-position}

The maximal Gaussian measure position provides both a lower covariance
bound and a Poincar\'e inequality at a fixed scale.

\begin{proposition}[Bobkov's Gaussian position]\label{prop:gaussian-position}
Every origin-symmetric convex body \(K\subset\R^n\) of volume one has a
maximal Gaussian measure position obtained by a map in \(\SL_n(\R)\).
In that position,
\begin{equation}\label{eq:alpha-covariance}
\Cov(\gamma_K)=\alpha\Id,
\qquad c_0\le\alpha\le1,
\end{equation}
where \(c_0>0\) is an absolute constant.
\end{proposition}

This is a consequence of \cite[Proposition~3.1 and Lemma~4.1]{Bob11}.
We include the proof of this proposition, since the covariance lower bound $c_0$ determines the
normalisation in our main proof.

\begin{proof}
The maximum exists because rotational invariance reduces the
maximisation to positive definite maps with determinant one. If such
maps leave every compact subset of this matrix space, their smallest
singular values tend to zero along a subsequence. Since the original
body is bounded, the images lie in slabs of widths tending to zero.
Their Gaussian measures therefore tend to zero. Continuity of the
Gaussian measure under nonsingular linear maps now gives a maximiser.

At a maximiser, vary the body by \(e^{tA}\), where \(A\) is symmetric
and \(\tr A=0\). Differentiation gives
\[
0=-\int_K\langle Ax,x\rangle\,\mathrm{d}\gam(x).
\]
Thus \(\int_K x\otimes x\,\mathrm{d}\gam(x)\) is a scalar matrix. Symmetry gives
zero mean, so \(\Cov(\gamma_K)=\alpha\Id\).

A consequence of Milman's \(M\)-ellipsoid theorem \cite{Mil86} is that
there exists a volume-one ellipsoid \(E\) with
\[
|K\cap E|\ge c_1^n.
\]
The equivalent intersection formulation is also described in
\cite[Corollary~2.3]{Bob11}. Send \(E\) by a determinant-one map to the
Euclidean ball \(B\) of volume one. Its squared radius is at most \(Cn\);
hence in this position
\[
\gam(K)\ge \gam(K\cap B)
\ge (2\pi)^{-n/2}e^{-Cn}|K\cap B|
\ge e^{-C_1n}.
\]
The maximising position has at least this Gaussian measure.

Let \(q=\gam(K)\) in the maximising position and \(X\sim\gamma_K\).
Its differential entropy is
\[
h(X)=\frac n2\log(2\pi)+\frac{n\alpha}{2}+\log q.
\]
The Gaussian maximises entropy at fixed covariance, so
\[
h(X)\le\frac n2\log(2\pi e\alpha).
\]
Consequently
\[
-C_1n\le\log q
\le\frac n2(1+\log\alpha-\alpha),
\]
which implies \(\alpha\ge e^{-2C_1-1}\).
Finally, the restricted Gaussian is \(1\)-strongly log-concave.
The Brascamp--Lieb Poincar\'e inequality implies
\(\Cov(X)\preceq\Id\), giving \(\alpha\le1\).
\end{proof}

For \(\kappa>0\), a density is called \(\kappa\)-strongly log-concave if
it has the form
\[
\rho(x)=e^{-\kappa|x|^2/2-W(x)}
\]
up to normalisation, with \(W:\R^n\to(-\infty,+\infty]\) convex.
Convex support is allowed. The Brascamp--Lieb inequality
\cite{BL76} gives
\begin{equation}\label{eq:strong-lc-poincare}
\CP(\rho\,\mathrm{d}x)\le\kappa^{-1},
\qquad
\Cov(\rho\,\mathrm{d}x)\preceq\kappa^{-1}\Id.
\end{equation}
For extended-valued \(W\), approximate \(W\) by finite smooth convex
potentials and pass to the variance inequality for compactly supported
smooth functions. The approximation gives \(C_P(\rho\,dx)\le\kappa^{-1}\).
Testing this inequality with \(\chi(x/R)(v\cdot x)\), where
\(\chi\in C_c^\infty(\mathbb R^n)\) equals one near the origin,
and letting \(R\to\infty\), we obtain
\[
v^{\mathsf{T}}\operatorname{Cov}(\rho\,dx)v
=\operatorname{Var}_{\rho\,dx}(v\cdot x)
\le\kappa^{-1}|v|^2.
\]
This holds for every \(v\in\mathbb R^n\), proving
\(\operatorname{Cov}(\rho\,dx)\preceq\kappa^{-1}I_n\).

We also record how Poincar\'e constants and Fisher information
change under convolution and linear maps.

\begin{lemma}\label{lem:convolution-facts}
Let \(X,Y\) be independent random vectors.
\begin{enumerate}[label=\textup{(\roman*)},leftmargin=2.2em]
\item If their Poincar\'e constants are finite, then
\[
\CP(X+Y)\le\CP(X)+\CP(Y).
\]
For a linear map \(A\),
\[
\CP(AX)\le\norm{A}_{\op}^2\CP(X).
\]
Consequently, normalised sums of independent copies of a law with
Poincar\'e constant at most \(P\) also have Poincar\'e constant at most \(P\).

\item If \(Y\) has finite ordinary Fisher information, then
\[
J(X+Y)\le J(Y).
\]
Here \(X\) may have an arbitrary law.
For an invertible matrix \(A\),
\[
J(AX)\le\norm{A^{-1}}_{\op}^2 J(X).
\]

\item Let \(U\sim\nu\) and \(G\sim\gam\) be independent, with
\(\nu\) isotropic. For \(t\ge0\), define
\[
\nu_t=\mathcal L(e^{-t}U+\sqrt{1-e^{-2t}}\,G).
\]
Then
\[
\CP(\nu_t)\le e^{-2t}\CP(\nu)+1-e^{-2t}.
\]
If \(I(\nu\mid\gam)<\infty\), then
\[
I(\nu_t\mid\gam)\le e^{-2t}I(\nu\mid\gam).
\]
\end{enumerate}
\end{lemma}

\begin{proof}
We begin by proving (i). For a smooth test function \(f\), conditioning on \(X\) gives
\[
\begin{aligned}
\operatorname{Var}(f(X+Y))
&=\mathbb{E}\bigl[\operatorname{Var}(f(X+Y)\mid X)\bigr]
  +\operatorname{Var}\bigl(\mathbb{E}[f(X+Y)\mid X]\bigr)\\
&\le \bigl(C_P(X)+C_P(Y)\bigr)
      \mathbb{E}|\nabla f(X+Y)|^2.
\end{aligned}
\]
Here we have used the Poincar\'e inequalities for \(Y\) and \(X\),
respectively, together with Jensen's inequality to bound the gradient
of the conditional mean. The bound for \(C_P(AX)\) follows by applying
the Poincar\'e inequality to \(f\circ A\) and using
\(\nabla(f\circ A)=A^{\mathsf{T}}(\nabla f)\circ A\).
This proves (i).

For (ii), write \(g\) for the density of \(Y\) and \(h\) for that of
\(X+Y\). The assumption \(J(Y)<\infty\) gives
\(\sqrt{g}\in H^1\), so
\[
\nabla g=2\sqrt{g}\,\nabla\sqrt{g}\in L^1.
\]
Differentiating the convolution defining \(h\) in the weak sense
therefore gives
\[
\nabla\log h(X+Y)
=\mathbb{E}[\nabla\log g(Y)\mid X+Y]
\]
on \(\{h>0\}\). Since conditional expectation is a contraction in
\(L^2\), we obtain \(J(X+Y)\le J(Y)\). A change of variables under
an invertible map \(A\) multiplies the score by \(A^{-\mathsf{T}}\),
giving
\[
J(AX)\le \|A^{-1}\|_{\mathrm{op}}^2J(X).
\]

The Poincar\'e bound in (iii) follows from (i), with the two summands
\(e^{-t}U\) and \(\sqrt{1-e^{-2t}}G\). The relative Fisher information
bound is the usual contraction estimate for the Ornstein--Uhlenbeck
(OU) semigroup; see \cite{BGL14}. Let \(P_t\) denote the OU semigroup.
The Fisher information estimate in (iii) is the standard contraction
bound for the Ornstein--Uhlenbeck (OU) semigroup:
\[
I(\nu_t\mid\gamma_n)\le e^{-2t}I(\nu\mid\gamma_n).
\]
It holds whenever \(I(\nu\mid\gamma_n)<\infty\); see \cite{BGL14}.
\end{proof}

\section{Second-order Stein kernels}
\label{sec:stein}

We construct a tensor satisfying the second-order Stein identity
for vector fields, following the spectral-gap argument in \cite{CFP19} and
its higher-order extension in \cite[Theorem~2.4]{Fathi21}. When \(\nu\) is even,
the construction gives \(T(-x)=-T(x)\), so
\(\int T\,d\nu=0\). The cross terms in the second-moment calculation
then vanish by independence, giving
\[
\|T_m\|_{L^2(\nu_m)}^2
\le m^{-2}\|T\|_{L^2(\nu)}^2,
\]
where \(T_m\) is the tensor constructed for the law \(\nu_m\) of
the normalised sum of \(m\) independent copies.

Throughout this section, the density of \(\nu\) is smooth and strictly
positive on \(\R^n\), \(\nu\) is even and isotropic, and
\(\CP(\nu)\le P\). Let \(U\sim\nu\).
In the application, these assumptions will be obtained
by Ornstein--Uhlenbeck regularisation of a compactly supported measure.

Let \(H^1(\nu)\) be the closure of \(C_c^\infty(\R^n)\) in the norm
\[
\left(\int |f|^2\,\mathrm{d}\nu+\int|\nabla f|^2\,\mathrm{d}\nu\right)^{1/2}.
\]
On each compact set, the density of \(\nu\) is bounded above and
below by positive constants. Thus, if \(f_k\to0\) and
\(\nabla f_k\to G\) in \(L^2(\nu)\), both convergences also hold
in \(L^2_{\mathrm{loc}}(dx)\). The local convergence identifies \(G\) as the weak gradient of
zero, so \(G=0\). This proves closability.

Let \(\chi\in C_c^\infty(\mathbb R^n)\) equal one on \(B_1\) and
vanish outside \(B_2\), and put \(\chi_R(x)=\chi(x/R)\). Then
\(\chi_R\to1\) in \(H^1(\nu)\). For the coordinate functions,
\[
\|x_j\chi_R-x_j\|_{L^2(\nu)}^2
+\|\nabla(x_j\chi_R)-e_j\|_{L^2(\nu)}^2
\le C\int_{|x|\ge R}(1+|x|^2)\,d\nu.
\]
The right-hand side tends to zero by the finite second moment of
\(\nu\), proving \(x_j\chi_R\to x_j\) in \(H^1(\nu)\).

On the mean-zero subspace, the Poincar\'e inequality gives
\[
\|f\|_{H^1(\nu)}^2
\le (1+P)\int |\nabla f|^2\,d\nu,
\]
so the quadratic form \(f\mapsto\int|\nabla f|^2\,d\nu\) is coercive.

\begin{lemma}\label{lem:second-stein}
There exists a third-order tensor \(T\in L^2(\nu)\), with
\[
\E T(U)=0,\qquad
\E\norm{T(U)}_{\HS}^2\le P(P-1)n,
\]
such that, for every \(F\in C_c^\infty(\R^n;\R^n)\),
\begin{equation}\label{eq:second-stein-identity}
\E[U\cdot F(U)-\Div F(U)]
=\E[T(U):D^2F(U)].
\end{equation}
The tensor \(T\) can be chosen odd.
\end{lemma}

\begin{proof}
For any mean-zero \(f\in H^1(\nu)\), Cauchy--Schwarz and the
Poincar\'e inequality give
\[
\left|\int x_i f\,d\nu\right|
\le \|x_i\|_{L^2(\nu)}\|f\|_{L^2(\nu)}
=\|f\|_{L^2(\nu)}
\le \sqrt{P}\,\|\nabla f\|_{L^2(\nu)},
\]
where isotropy gives \(\|x_i\|_{L^2(\nu)}=1\). Thus the functional
\(f\mapsto\int x_i f\,d\nu\) is bounded in the gradient norm.
Lax--Milgram yields a unique mean-zero \(u_i\in H^1(\nu)\) satisfying
\begin{equation}\label{eq:first-poisson}
\int \nabla u_i\cdot\nabla f\,d\nu
=\int x_i f\,d\nu
\end{equation}
for every mean-zero \(f\in H^1(\nu)\). Since \(\int x_i\,d\nu=0\), the identity holds for constants as
well as mean-zero functions, and hence for every \(f\in H^1(\nu)\).
Taking \(f=u_i\), Cauchy--Schwarz and the Poincar\'e inequality give
\[
\int|\nabla u_i|^2\,d\nu
=\int x_i u_i\,d\nu
\le \|u_i\|_{L^2(\nu)}
\le \sqrt{P}\left(\int|\nabla u_i|^2\,d\nu\right)^{1/2},
\]
and therefore \(\int|\nabla u_i|^2\,d\nu\le P\).
Define \(\tau_{ij}=\partial_j u_i\). Testing
 with \(x_j\) gives
\(\E\tau_{ij}=\delta_{ij}\). Therefore
\begin{equation}\label{eq:first-stein-size}
\E\norm{\tau-\Id}_{\HS}^2
=\sum_i\int|\nabla u_i|^2\,\mathrm{d}\nu-n
\le(P-1)n.
\end{equation}
The equation also gives the first-order Stein identity
\[
\E[U\cdot F(U)]=\E[\tau(U):DF(U)].
\]

The measure is even. Reflection of
\eqref{eq:first-poisson} and uniqueness show that \(u_i\) is odd.
Thus \(\tau\) is even. For each pair \(i,j\), solve the mean-zero weak
equation
\begin{equation}\label{eq:second-poisson}
\int\nabla v_{ij}\cdot\nabla f\,\mathrm{d}\nu
=\int(\tau_{ij}-\delta_{ij})f\,\mathrm{d}\nu.
\end{equation}
Because \(\tau_{ij}-\delta_{ij}\) has mean zero and belongs to
\(L^2(\nu)\), Lax--Milgram gives a unique mean-zero solution
\(v_{ij}\) of \eqref{eq:second-poisson}. Testing with \(v_{ij}\)
and using the Poincar\'e inequality gives
\[
\int|\nabla v_{ij}|^2\,\mathrm{d}\nu
\le P\int|\tau_{ij}-\delta_{ij}|^2\,\mathrm{d}\nu.
\]

Since \(\nu\) and \(\tau\) are even, \(x\mapsto v_{ij}(-x)\)
is another mean-zero solution of \eqref{eq:second-poisson}.
Uniqueness therefore implies that \(v_{ij}\) is even. The tensor
\(T\) defined by
\[
T_{ijk}=\partial_k v_{ij}
\]
is consequently odd and has mean zero. Its norm satisfies
\[
\begin{aligned}
\E\norm T_{\HS}^2
&=\sum_{i,j}\int|\nabla v_{ij}|^2\,\mathrm{d}\nu\\
&\le P\,\E\norm{\tau-I_n}_{\HS}^2
\le P(P-1)n,
\end{aligned}
\]
where the last inequality is \eqref{eq:first-stein-size}.

For a smooth, compactly supported vector field \(F\), taking
\(f=\partial_jF_i\) in \eqref{eq:second-poisson} and summing over
\(i,j\) gives
\[
\int T:D^2F\,\mathrm{d}\nu
=\int(\tau-I_n):DF\,\mathrm{d}\nu
=\int\bigl(x\cdot F-\operatorname{div}F\bigr)\,\mathrm{d}\nu.
\]
The last equality is the first-order Stein identity, which proves
\eqref{eq:second-stein-identity}.
\end{proof}

The following construction is the second-order case of
\cite[Lemma~4.1 and Corollary~4.2]{Fathi21}.

\begin{lemma}[Normalised sums]\label{lem:stein-sums}
Let \(U_1,\ldots,U_m\) be independent with law \(\nu\), and define
\[
W_m=m^{-1/2}\sum_{\ell=1}^m U_\ell.
\]
Then the tensor
\begin{equation}\label{eq:sum-tensor}
T_m(W_m)=m^{-3/2}
\E\left[\sum_{\ell=1}^m T(U_\ell)\,\middle|\,W_m\right]
\end{equation}
satisfies the second-order Stein identity for \(W_m\), and
\begin{equation}\label{eq:sum-tensor-bound}
\E\norm{T_m(W_m)}_{\HS}^2
\le\frac{P(P-1)n}{m^2}.
\end{equation}
\end{lemma}

\begin{proof}
For each \(\ell\), applying \eqref{eq:second-stein-identity}
conditionally on the other summands, with test field
\[
x\longmapsto F\left(\frac{x+\sum_{j\ne\ell}U_j}{\sqrt m}\right),
\]
gives
\[
\E\left[U_\ell\cdot F(W_m)
-\frac{1}{\sqrt m}\operatorname{div}F(W_m)\right]
=\frac{1}{m}\E\left[T(U_\ell):D^2F(W_m)\right].
\]
Multiplying by \(m^{-1/2}\) and summing over \(\ell\) gives the
coefficient \(m^{-3/2}\) in \eqref{eq:sum-tensor}..

Conditional Jensen's inequality and independence give
\[
\E\norm{T_m(W_m)}_{\HS}^2
\le m^{-3}\E\left\|\sum_{\ell=1}^m T(U_\ell)\right\|_{\HS}^2.
\]
Since \(T\) is centred, all cross terms vanish, and the last expression
is \(m^{-2}\E\norm T_{\HS}^2\). Apply
Lemma~\ref{lem:second-stein}.
\end{proof}

We next bound the entropy in terms of Fisher information and the
\(L^2\) norm of the second-order Stein tensor. We adapt the proof
of \cite[Theorem~3.1 and Lemma~3.2]{Fathi21} to the vector-field identity above,
using the semigroup method described in \cite{LNP15}.

\begin{lemma}[Entropy from a second-order tensor]\label{lem:hsi}
Let \(U\) have an isotropic law \(\nu\) with finite relative Fisher
information and admit a tensor \(T\in L^2(\nu)\) satisfying
\eqref{eq:second-stein-identity}. Then
\begin{equation}\label{eq:hsi}
H(\nu\mid\gam)
\le\sqrt2\,\norm T_{L^2(\nu)}\sqrt{I(\nu\mid\gam)}.
\end{equation}
\end{lemma}

\begin{proof}
For \(t>0\), let \(a=e^{-t}\), \(b=\sqrt{1-e^{-2t}}\), and
\(U_t=aU+bG\), where
\(G\sim\gam\) is independent. Denote the relative score of its law by
\[
r_t=\nabla\log\frac{d\nu_t}{d\gam}.
\]
We first work against a vector field \(F\in C_c^\infty(\R^n;\R^n)\).
Integration by parts for the relative score gives
\[
\E[r_t(U_t)\cdot F(U_t)]
=\E[U_t\cdot F(U_t)-\Div F(U_t)].
\]
Gaussian integration by parts in the term containing \(bG\), followed
by \eqref{eq:second-stein-identity} in \(U\), gives
\begin{align*}
\E[U_t\cdot F(U_t)-\Div F(U_t)]
&=a\,\E[U\cdot F(aU+bG)]-a^2\E[\Div F(aU+bG)]\\
&=a^3\E[T(U):D^2F(aU+bG)].
\end{align*}
The Stein identity applies to \(x\mapsto F(ax+bG)\), which is
smooth and compactly supported for each fixed \(G\). Integrating
twice more by parts with respect to \(G\) gives
\[
\E[r_t(U_t)\cdot F(U_t)]
=\frac{a^3}{b^2}
\E\bigl[(T(U):(G\otimes G-\Id))\cdot F(U_t)\bigr].
\]
The derivatives of \(F\) are bounded, and \(T\in L^2(\nu)\), which
justifies these integrations. Equivalently, the relative score is
identified weakly as
\begin{equation}\label{eq:ou-score}
r_t(U_t)
=\frac{a^3}{b^2}
\E[T(U):(G\otimes G-\Id)\mid U_t].
\end{equation}

For a fixed tensor \(T\), the Gaussian fourth-moment formula gives
\begin{equation}\label{eq:hermite-contraction}
\E_G|T:(G\otimes G-\Id)|^2
=\sum_{i,j,k}\bigl(T_{ijk}^2+T_{ijk}T_{ikj}\bigr)
\le2\norm T_{\HS}^2.
\end{equation}
Thus \eqref{eq:ou-score} and conditional Jensen's inequality imply
\begin{equation}\label{eq:ou-second-bound}
I(\nu_t\mid\gam)
\le\frac{2e^{-6t}}{(1-e^{-2t})^2}\norm T_{L^2(\nu)}^2.
\end{equation}
Together with the relative Fisher-information contraction under the OU
semigroup, this gives
\[
I(\nu_t\mid\gam)
\le
\min\left\{
e^{-2t}I(\nu\mid\gam),
\frac{2e^{-6t}}{(1-e^{-2t})^2}\norm T_{L^2(\nu)}^2
\right\}.
\]

The entropy dissipation identity for the OU semigroup is
\[
H(\nu\mid\gam)=\int_0^\infty I(\nu_t\mid\gam)\,\mathrm{d}t;
\]
see, for example, \cite{BGL14}. Finite information ensures finite
entropy by the Gaussian logarithmic Sobolev inequality. The same
inequality gives
\(H(\nu_t\mid\gam)\le I(\nu_t\mid\gam)/2\to0\) as \(t\to\infty\).
At \(t\downarrow0\), lower semicontinuity of entropy and its contraction
under the semigroup give \(H(\nu_t\mid\gam)\to H(\nu\mid\gam)\).
Thus both endpoints in the dissipation identity are justified. Set
\(I_*=I(\nu\mid\gam)\) and \(s=\norm T_{L^2(\nu)}\).
If \(I_*=0\), the Gaussian logarithmic Sobolev inequality gives
\(H(\nu\mid\gam)=0\).
If \(s=0\), \eqref{eq:ou-second-bound} and entropy dissipation give
\(H(\nu\mid\gam)=0\). We may therefore assume \(I_*,s>0\). With
\(z=1-e^{-2t}\),
\[
H(\nu\mid\gam)
\le\frac12\int_0^1
 \min\left\{I_*,\frac{2s^2(1-z)^2}{z^2}\right\}\,\mathrm{d}z
\le\frac12\int_0^1
 \min\left\{I_*,\frac{2s^2}{z^2}\right\}\,\mathrm{d}z.
\]
If \(\sqrt2s/\sqrt{I_*}\le1\), split the integral at this value. The
result is at most \(\sqrt2s\sqrt{I_*}\). If
\(\sqrt2s/\sqrt{I_*}>1\), the integral is at most \(I_*/2\), which
satisfies the same upper bound. This proves
\eqref{eq:hsi}.
\end{proof}

\section{Entropy decay and Gaussian measure}
\label{sec:entropy}

The Stein estimate supplies one power of \(m^{-1}\) in the entropy of a
normalised sum. Combining this with Bizeul's information-jump inequality gives
entropy decay faster than \(m^{-1}\), as required for the Brownian
bridge argument.

\begin{proposition}[Bizeul's information jump]\label{prop:fisher-jump}
Let \(U,U'\) be independent and identically distributed, with an
isotropic log-concave law \(\nu\) and
\(\CP(\nu)\le P\). If \(I(\nu\mid\gam)<\infty\), then
\begin{equation}\label{eq:fisher-jump}
I\left(\mathcal L\left(\frac{U+U'}{\sqrt2}\right)\middle|\gam\right)
\le \left(1-\frac1{8P}\right)I(\nu\mid\gam).
\end{equation}
\end{proposition}

\begin{proof}
This is \cite[Theorem~2]{Biz23}, with equal input laws and mixing
parameter \(1/2\). The deficit in that theorem is expressed in ordinary
Fisher information, while its lower bound uses relative Fisher
information. The output and input laws are isotropic, so
\eqref{eq:relative-ordinary} shows that their ordinary and relative
information deficits coincide. The resulting lower bound on the
deficit is \(I(\nu\mid\gam)/(8P)\).

For a nonsmooth density, apply the estimate first to
\(\nu_\eps=\mathcal L(e^{-\eps}U+\sqrt{1-e^{-2\eps}}\,G)\).
This law is smooth, positive, isotropic and log-concave, and
Lemma~\ref{lem:convolution-facts} gives
\(\CP(\nu_\eps)\le P\) and
\(I(\nu_\eps\mid\gam)\le I(\nu\mid\gam)\).
Taking a normalised sum of two independent copies commutes with this
regularisation. Lower semicontinuity of relative Fisher information at
the output as \(\eps\downarrow0\) proves \eqref{eq:fisher-jump}.
\end{proof}

Let \(I_0=I(\nu\mid\gam)<\infty\), and let \(\nu_m\) be the law of
the normalised sum of \(m\) independent copies with law \(\nu\).
Log-concavity follows
from the Pr\'ekopa--Leindler theorem, isotropy is preserved, and
Lemma~\ref{lem:convolution-facts} gives \(\CP(\nu_m)\le P\).
Consequently, for dyadic \(m\),
\begin{equation}\label{eq:dyadic-fisher}
I(\nu_m\mid\gam)\le m^{-\delta}I_0,
\qquad
\delta=-\log_2\left(1-\frac1{8P}\right)>0.
\end{equation}

\begin{proposition}[Entropy of symmetric normalised sums]
\label{prop:entropy-clt}
Suppose \(\nu\) is an even isotropic log-concave probability supported
on a convex body, with \(\CP(\nu)\le P\), \(P\ge2\), and finite
\(I_0=I(\nu\mid\gam)\). Then, for dyadic \(m\),
\begin{equation}\label{eq:entropy-clt}
H(\nu_m\mid\gam)
\le \sqrt{2P(P-1)n I_0}\,m^{-1-\delta/2},
\end{equation}
where \(\delta\) is defined in \eqref{eq:dyadic-fisher}.
\end{proposition}

\begin{proof}
Assume first that \(\nu\) has a smooth positive density on \(\R^n\),
is even, isotropic and log-concave, and satisfies the stated Poincar\'e
and information bounds.
Lemmas~\ref{lem:stein-sums} and \ref{lem:hsi}, together with
\eqref{eq:dyadic-fisher}, give
\[
H(\nu_m\mid\gam)
\le \sqrt2\,\norm{T_m}_{L^2(\nu_m)}
       \sqrt{I(\nu_m\mid\gam)}
\le \sqrt{2P(P-1)n I_0}\,m^{-1-\delta/2}.
\]

For the compactly supported law in the statement, set
\[
\nu_\eps
=\mathcal L(e^{-\eps}U+\sqrt{1-e^{-2\eps}}\,G).
\]
The density of \(\nu_\eps\) is smooth, strictly positive, even, and
log-concave. It remains isotropic, and
Lemma~\ref{lem:convolution-facts} gives
\[
\CP(\nu_\eps)\le e^{-2\eps}P+1-e^{-2\eps}\le P,
\qquad
I(\nu_\eps\mid\gam)\le I_0.
\]
Since \(\nu\) has compact support, \(\nu_\varepsilon\) has Gaussian
tails and hence finite moments of every order. Hence, this justifies the cutoff
approximations and integrations by parts used in Section~3.
A normalised sum of independent standard Gaussian vectors is again
standard Gaussian, so for every fixed dyadic \(m\),
\[
(\nu_\varepsilon)_m=(\nu_m)_\varepsilon.
\]
The preceding entropy estimate therefore gives a bound for
\(H((\nu_m)_\varepsilon\mid\gamma_n)\). This bound is independent
of \(\varepsilon\), since
\(C_P(\nu_\varepsilon)\le P\) and
\(I(\nu_\varepsilon\mid\gamma_n)\le I_0\).
As \(\eps\downarrow0\), weak convergence and lower semicontinuity of
relative entropy therefore yield
\[
H(\nu_m\mid\gam)
\le\liminf_{\eps\downarrow0}H((\nu_\eps)_m\mid\gam)
\le\sqrt{2P(P-1)n I_0}\,m^{-1-\delta/2}.
\]
\end{proof}

\begin{lemma}[Entropy and support]\label{lem:support-entropy}
If a probability measure \(\eta\) is supported on a measurable set \(A\)
of positive Gaussian measure, then
\[
-\log\gam(A)\le H(\eta\mid\gam).
\]
\end{lemma}

\begin{proof}
If the relative entropy is infinite there is nothing to prove.
Otherwise \(\eta\ll\gam(\cdot\mid A)\), and
\[
H(\eta\mid\gam)
=H(\eta\mid\gam(\cdot\mid A))-\log\gam(A)
\ge-\log\gam(A).
\]
\end{proof}

\begin{proposition}\label{prop:gaussian-profile}
Under the assumptions of Theorem~\ref{thm:information-spectrum}, let
\(p=2+\delta>2\). There is
\[
A\le C(P)\sqrt{n I_0}
\]
such that
\begin{equation}\label{eq:gaussian-profile}
-\log\gam(tD)\le A t^{-p}\qquad(t\ge1).
\end{equation}
\end{proposition}

\begin{proof}
Convexity gives
\(\supp(\nu_m)\subset m^{-1/2}(D+\cdots+D)=\sqrt mD\).
Apply Lemma~\ref{lem:support-entropy} and
Proposition~\ref{prop:entropy-clt} to obtain, for dyadic \(m\),
\[
-\log\gam(\sqrt mD)
\le \sqrt{2P(P-1)n I_0}\,m^{-p/2}.
\]
For a given \(t\ge1\), choose the largest dyadic integer \(m\le t^2\).
Then \(t^2/2<m\le t^2\). Monotonicity under dilation gives
\[
-\log\gam(tD)
\le 2^{p/2}\sqrt{2P(P-1)n I_0}\,t^{-p}.
\]
This proves the assertion.
\end{proof}

\section{Brownian bridges and the first eigenvalue}
\label{sec:bridge}

The preceding section gives a bound on Gaussian measure under dilation.
We now convert that bound into a lower bound for Brownian bridge
survival, and hence an upper bound for the first Dirichlet eigenvalue.

\begin{lemma}
\label{lem:bridge}
Let \(D\subset\R^n\) be an origin-symmetric convex body. Suppose that
\(p>2\) and \(A\ge0\) satisfy
\[
-\log\gam(tD)\le A t^{-p}\qquad(t\ge1).
\]
Then
\begin{equation}\label{eq:bridge-eigenvalue}
\lambda_1(D)\le C(p)(A+n).
\end{equation}
\end{lemma}

\begin{proof}
Let \(\beta\) be the standard \(n\)-dimensional Brownian bridge from
zero to zero on \([0,1]\), associated with Brownian motion of generator
\(\Delta/2\). Its L\'evy--Ciesielski expansion is
\begin{equation}\label{eq:bridge-expansion}
\beta(t)=\sum_{j=0}^\infty 2^{-j/2-1}
 \sum_{k=0}^{2^j-1}G_{jk}h_{jk}(t),
\end{equation}
where the \(G_{jk}\) are independent standard Gaussian vectors.
The function \(h_{jk}\) is supported on
\([k2^{-j},(k+1)2^{-j}]\), is linear on each half of that interval,
vanishes at its endpoints, and equals one at the midpoint.
At each level \(j\) and time \(t\), at most one tent is nonzero.

The coefficients are the normalised midpoint residuals
\[
2^{j/2+1}\left[
 \beta\bigl((2k+1)2^{-j-1}\bigr)
 -\frac12\left\{\beta(k2^{-j})+\beta((k+1)2^{-j})\right\}
\right].
\]
Conditional midpoint variances are \(2^{-j-2}\), and Gaussian
conditioning gives independence of the residuals. This yields
\eqref{eq:bridge-expansion}; the standard L\'evy construction is discussed
in \cite[Section~1.1.2]{MP10}.

Choose a fixed number \(a\) satisfying
\begin{equation}\label{eq:choice-a}
\frac1p<a<\frac12.
\end{equation}
Set \(b_j=2^{aj}\) and consider the coefficient event
\[
E=\bigcap_{j\ge0}\ \bigcap_{k=0}^{2^j-1}\{G_{jk}\in b_jD\}.
\]
By independence,
\begin{align}
-\log\Prob(E)
&=\sum_{j=0}^\infty 2^j[-\log\gam(b_jD)]\notag\\
&\le A\sum_{j=0}^\infty2^{(1-ap)j}
=\frac{A}{1-2^{1-ap}}.
\label{eq:bridge-probability}
\end{align}
For finite intersections, independence gives the corresponding finite
product. The finite intersections decrease to \(E\), so their probabilities
converge to \(\mathbb{P}(E)\), giving the product formula. Since
\(ap>1\), the negative logarithms of the factors are summable, and
the product is strictly positive.

Let \(\norm{\cdot}_D\) be the norm with unit ball \(D\).
On \(E\), \eqref{eq:bridge-expansion} gives the deterministic uniform
bound
\begin{equation}\label{eq:bridge-radius}
\sup_{t\in[0,1]}\norm{\beta(t)}_D
\le\frac12\sum_{j=0}^\infty2^{(a-1/2)j}
=R,\qquad
R=\frac1{2(1-2^{a-1/2})}.
\end{equation}
This same convergent series controls uniform convergence of the
expansion on \(E\). The expansion agrees almost surely with the
continuous bridge.

Let \(\Omega=2R \Int D\).
Since \(0\in \Int D\), the compact set \(RD\) is contained in \(\Omega\)
with a positive distance from its boundary. Thus \(E\) implies
survival of the bridge inside the open domain. If
\[
q=\Prob\{\beta(t)\in\Omega\text{ for every }t\in[0,1]\},
\]
then
\begin{equation}\label{eq:q-lower}
q\ge\exp\left(-\frac{A}{1-2^{1-ap}}\right).
\end{equation}

Let \(p_t^\Omega\) be the Dirichlet heat kernel for the generator
\(\Delta/2\). The killed Brownian transition representation and bridge
conditioning give
\begin{equation}\label{eq:bridge-kernel}
p_1^\Omega(0,0)=q\,p_1^{\R^n}(0,0).
\end{equation}
See \cite[Section~7.4]{MP10} for the killed-process representation.

Let \(\lambda_\ell(\Omega)\) and \(\phi_\ell\) be Dirichlet
eigenvalues of \(-\Delta\) and an orthonormal eigenbasis. Positivity of
the squared coefficients in the diagonal spectral expansion yields
\begin{align}
p_1^\Omega(0,0)
&=\sum_\ell e^{-\lambda_\ell(\Omega)/2}\phi_\ell(0)^2\notag\\
&\le e^{-\lambda_1(\Omega)/4}
  \sum_\ell e^{-\lambda_\ell(\Omega)/4}\phi_\ell(0)^2\notag\\
&=e^{-\lambda_1(\Omega)/4}p_{1/2}^\Omega(0,0)\notag\\
&\le e^{-\lambda_1(\Omega)/4}p_{1/2}^{\R^n}(0,0).
\label{eq:heat-kernel-bound}
\end{align}
The heat kernel and eigenfunctions are smooth in the interior of
\(\Omega\), so these diagonal evaluations are well defined. The last
inequality is domination by the free heat kernel.
For Brownian motion of generator \(\Delta/2\),
\(p_t^{\R^n}(0,0)=(2\pi t)^{-n/2}\), so
\eqref{eq:bridge-kernel}--\eqref{eq:heat-kernel-bound} give
\[
q\le 2^{n/2}e^{-\lambda_1(\Omega)/4}.
\]
It follows that
\[
\lambda_1(\Omega)
\le4\log(1/q)+2n\log2
\le\frac{4A}{1-2^{1-ap}}+2n\log2.
\]
Finally,
\[
\lambda_1(D)=(2R)^2\lambda_1(\Omega).
\]
Choosing, for example,
\(a=\tfrac12(\tfrac1p+\tfrac12)\) makes every resulting constant
depend only on \(p\), which proves \eqref{eq:bridge-eigenvalue}.
\end{proof}

\begin{proof}[Proof of Theorem~\ref{thm:information-spectrum}]
Proposition~\ref{prop:gaussian-profile} gives
\[
-\log\gam(tD)\le C(P)\sqrt{nI_0}\,t^{-p},
\qquad
p=2-\log_2\left(1-\frac1{8P}\right)>2.
\]
Apply Lemma~\ref{lem:bridge}. Both the exponent \(p\) and the
constants depend only on \(P\), yielding
\eqref{eq:information-spectrum}.
\end{proof}

\begin{remark}\label{rem:critical-exponent}
The two series in \eqref{eq:bridge-probability} and
\eqref{eq:bridge-radius} converge simultaneously when
\(1/p<a<1/2\), which is possible exactly when \(p>2\).
The centred Stein tensor contributes a factor \(m^{-1}\) to the entropy
estimate. The information jump adds \(m^{-\delta/2}\), giving the
Gaussian dilation exponent \(2+\delta\).
\end{remark}

\section{Proof of the affine reverse Faber--Krahn inequality}
\label{sec:main-proof}

Theorem~\ref{thm:information-spectrum} reduces the spectral estimate to
the construction of a measure with controlled Poincar\'e constant and
Fisher information. We use the ground state to make the density vanish at the boundary,
then convolve with the restricted Gaussian to obtain a uniform
lower bound on the covariance.

Let \(K\) be origin-symmetric, of volume one, in maximal Gaussian
measure position. Let
\[
\lambda=\lambda_1(K),\qquad
\Cov(\gamma_K)=\alpha\Id,\qquad c_0\le\alpha\le1.
\]
Let \(\psi\in H_0^1(\Int K)\) be a positive first eigenfunction.
It is unique up to a positive scalar. Since \(K=-K\), it can be
chosen even. The Brascamp--Lieb log-concavity theorem for the first
Dirichlet eigenfunction on a convex domain \cite{BL76} gives
log-concavity of \(\psi\). The extension by zero is understood in this
log-concavity statement.

\begin{lemma}[Gaussian-weighted ground state]\label{lem:groundstate}
Define
\[
v(x)=\psi(x)e^{-|x|^2/4}\quad(x\in \Int K),
\]
and extend \(v\) by zero. Let \(Y\) have probability density
\[
q_Y(x)=\frac{v(x)^2}{\int_{\R^n}v^2}.
\]
Then \(Y\) is even and \(1\)-strongly log-concave, and
\[
\CP(Y)\le1,\qquad
\Cov(Y)\preceq\Id,\qquad
J(Y)\le4\lambda+n.
\]
More precisely,
\begin{equation}\label{eq:groundstate-exact}
J(Y)=4\lambda+\E|Y|^2.
\end{equation}
\end{lemma}

\begin{proof}
The density is proportional to
\(\psi^2e^{-|x|^2/2}\), so log-concavity of \(\psi\) gives
\(1\)-strong log-concavity. Symmetry gives zero mean.
Equation~\eqref{eq:strong-lc-poincare} supplies the Poincar\'e and
covariance bounds.

To compute the information, test the weak equation
\(-\Delta\psi=\lambda\psi\) with
\(\psi e^{-|x|^2/2}\). This is an admissible member of
\(H_0^1(\Int K)\). With \(w=e^{-|x|^2/2}\), this gives
\[
\int_K w|\nabla\psi|^2
-\int_K w\psi\,x\cdot\nabla\psi
=\lambda\int_K w\psi^2.
\]
On the other hand,
\[
|\nabla v|^2
=w|\nabla\psi|^2-w\psi\,x\cdot\nabla\psi
+\frac14|x|^2w\psi^2.
\]
Thus
\begin{equation}\label{eq:groundstate-energy}
\int|\nabla v|^2
=\lambda\int v^2+\frac14\int|x|^2v^2.
\end{equation}
Multiplication by the smooth Gaussian weight preserves \(H_0^1\), and
zero extension places \(v\) in \(H^1(\R^n)\). Therefore the definition
\eqref{eq:ordinary-fisher} applies and gives
\[
J(Y)=4\frac{\int|\nabla v|^2}{\int v^2}
=4\lambda+\E|Y|^2.
\]
The covariance bound and evenness imply \(\E|Y|^2\le n\).
\end{proof}

Independently of \(Y\), let \(X\sim\gamma_K\), and define
\begin{equation}\label{eq:constructed-measure}
\Sigma=\Cov(X+Y)=\alpha\Id+\Cov(Y),\qquad
Z=\Sigma^{-1/2}(X+Y),\qquad
D=2\Sigma^{-1/2}K.
\end{equation}
By Proposition~\ref{prop:gaussian-position} and
Lemma~\ref{lem:groundstate},
\begin{equation}\label{eq:sigma-bounds}
\alpha\Id\preceq\Sigma\preceq2\Id.
\end{equation}

\begin{lemma}\label{lem:constructed-bounds}
The law \(\nu=\mathcal L(Z)\) is even, isotropic, log-concave and
supported on \(D\). It satisfies
\begin{equation}\label{eq:constructed-poincare}
\CP(\nu)\le\frac2\alpha\le\frac2{c_0},
\end{equation}
and
\begin{equation}\label{eq:constructed-fisher}
I(\nu\mid\gam)\le8\lambda+n.
\end{equation}
\end{lemma}

\begin{proof}
Evenness, log-concavity and support follow from independence,
convolution of log-concave densities, and \(K+K=2K\).
The definition of \(\Sigma\) gives isotropy.
The law of \(X\) is \(1\)-strongly log-concave, so
\(\CP(X)\le1\). Lemma~\ref{lem:groundstate} gives
\(\CP(Y)\le1\). Hence Lemma~\ref{lem:convolution-facts} implies
\[
\CP(Z)\le\norm{\Sigma^{-1/2}}_{\op}^2\CP(X+Y)
\le\frac2\alpha.
\]

The same lemma, now for Fisher information, gives
\[
J(Z)\le\norm{\Sigma}_{\op}J(X+Y)
\le2J(Y)\le8\lambda+2n.
\]
Here only \(Y\) is required to have finite Fisher information; the
density of \(X\) may have a jump at \(\partial K\).
Since \(Z\) is isotropic, \eqref{eq:relative-ordinary} proves
\eqref{eq:constructed-fisher}.
\end{proof}

\begin{lemma}\label{lem:linear-rayleigh}
For any invertible \(A:\R^n\to\R^n\) and any bounded domain \(U\),
\[
\lambda_1(U)\le\norm A_{\op}^2\lambda_1(AU).
\]
\end{lemma}

\begin{proof}
For \(f\in H_0^1(AU)\), use \(f\circ A\in H_0^1(U)\) as a test
function. Change of variables gives
\[
\frac{\int_U|\nabla(f\circ A)|^2}{\int_U|f\circ A|^2}
=
\frac{\int_{AU}|A^{\mathsf T}\nabla f|^2}{\int_{AU}|f|^2}
\le\norm A_{\op}^2
\frac{\int_{AU}|\nabla f|^2}{\int_{AU}|f|^2}.
\]
Take the infimum over nonzero \(f\).
\end{proof}

\begin{proof}[Proof of Theorem~\ref{thm:gaussian-main}]
Use the measure in \eqref{eq:constructed-measure}.
The fixed number \(P_0=2/c_0\) is at least two and is universal.
Theorem~\ref{thm:information-spectrum} and
Lemma~\ref{lem:constructed-bounds} give
\[
\lambda_1(D)
\le C(P_0)\left(n+\sqrt{n(8\lambda+n)}\right).
\]
For \(A=2\Sigma^{-1/2}\), equation~\eqref{eq:sigma-bounds} gives
\(\norm A_{\op}^2\le4/\alpha\le4/c_0\).
Since \(D=AK\), Lemma~\ref{lem:linear-rayleigh} yields
\[
\lambda
\le C\left(n+\sqrt{n(8\lambda+n)}\right)
\le C'\left(n+\sqrt{n\lambda}\right),
\]
with absolute constants. By Young's inequality,
\[
C'\sqrt{n\lambda}
\le\frac{\lambda}{2}+\frac{(C')^2n}{2}.
\]
Absorbing \(\lambda/2\) proves
\[
\lambda_1(K)\le\bigl(2C'+(C')^2\bigr)n.
\]
\end{proof}

\subsection{Reduction to symmetric bodies}
\label{sec:general-case}

To pass from the symmetric case to a general convex body \(K\),
we use a centrally symmetric convex body \(L\subset K\) with
\(|L|\ge2^{-n}|K|\). The averaging argument below shows that such
a body exists; see also \cite[Section~1.2]{HSTV22}.

\begin{proof}[Proof of Theorem~\ref{thm:affine-main}]
Write \(V=|K|\). If \(K=-K\), choose \(S\in\SL_n(\R)\) so that the
volume-one body \(V^{-1/n}SK\) is in maximal Gaussian measure
position. Theorem~\ref{thm:gaussian-main} and the scaling of the
Dirichlet eigenvalue give
\[
\lambda_1(SK)V^{2/n}
=\lambda_1(V^{-1/n}SK)
\le Cn.
\]

For a general convex body \(K\), Fubini's theorem gives
\[
\int_{\R^n}|K\cap(x-K)|\,\mathrm{d}x=V^2.
\]
Since the integrand is supported on \(K+K=2K\), a set of volume
\(2^nV\), there is a point \(x_0\in\R^n\) for which
\[
|K\cap(x_0-K)|\ge2^{-n}V.
\]
The intersection \(L=K\cap(x_0-K)\) is a convex body contained in
\(K\) and symmetric about \(x_0/2\). Hence we can apply the
symmetric case to \(L-x_0/2\). By translation invariance, this gives
a map \(S\in\SL_n(\R)\) satisfying
\[
\lambda_1(SL)|L|^{2/n}\le Cn.
\]
Because \(SL\subset SK\), domain monotonicity now yields
\[
\lambda_1(SK)V^{2/n}
\le \lambda_1(SL)V^{2/n}
\le Cn\left(\frac{V}{|L|}\right)^{2/n}
\le4Cn,
\]
which proves the result after enlarging the universal constant.
\end{proof}

\paragraph{\textbf{AI declaration}}
The authors interacted with and guided ChatGPT 6 whilst developing the proofs and then
wrote the paper. After we wrote the paper, ChatGPT 6 was used for
minor language editing of the penultimate draft in this stage only. The authors take full responsibility for the accuracy and integrity of all
content in this paper

\end{document}